\documentclass[12pt,a4paper]{amsart}
\usepackage{amsmath,latexsym,epsfig,graphics,amsfonts,amssymb,amsthm}

\newtheorem{theorem}{Theorem}
\newtheorem{lemma}{Lemma}[section]
\newtheorem{proposition}[lemma]{Proposition}

\theoremstyle{definition}

\newtheorem{example}[lemma]{Example}

\theoremstyle{remark}

\theoremstyle{plain}

\numberwithin{equation}{section}
\numberwithin{figure}{section}

\renewcommand{\leq}{\leqslant}
\renewcommand{\geq}{\geqslant}
\newcommand{\la}{\ensuremath{\langle}}
\newcommand{\ra}{\ensuremath{\rangle}}
\newcommand{\RN}{\ensuremath{\mathbb{R}^N_\infty}}

\begin{document}

\title{Rigidity of configurations of balls and points in the $N$-sphere}


\author{Edward Crane}
\address{Heilbronn Institute for Mathematical Research \\
Royal Fort Annexe\\
University Walk \\
Bristol \\
BS8 1TW\\
United Kingdom}
\email{edward.crane@gmail.com}


\author{Ian Short}
\address{Department of Pure Mathematics and Mathematical Statistics\\
Centre for Mathematical Sciences \\
Wilberforce Road \\
Cambridge\\
CB3 0WB\\
United Kingdom}
\email{ims25@cam.ac.uk}

\date{\today}

\subjclass[2000]{Primary: 51B10; Secondary: 30C20}
\keywords{Automorphism, ball, circular region, conformal, hyperbolic geometry, hyperboloid model, rigidity, sphere, symmetry.}

\begin{abstract}
We answer two questions of Beardon and Minda that arose from their study of the conformal symmetries of circular regions in the complex plane. We show that a configuration of closed balls in the $N$-sphere is determined up to M\"{o}bius transformations by the signed inversive distances between pairs of its elements, except when the boundaries of the balls have a point in common, and that a configuration of points in the $N$-sphere is determined by the absolute cross-ratios of $4$-tuples of its elements. The proofs use the hyperboloid model of hyperbolic $(N+1)$-space.
\end{abstract}

\maketitle

\section{Introduction}\label{S: introduction}

Let $\mathbb{R}^N_\infty$ denote the one point extension of $\mathbb{R}^N$. A \emph{sphere} in  $\mathbb{R}^N_\infty$ refers to either an $(N-1)$-dimensional  Euclidean sphere in $\mathbb{R}^N$, or an $(N-1)$-dimensional Euclidean plane  in $\mathbb{R}^N$ with the point $\infty$ attached. An \emph{inversion} in a sphere refers to either a Euclidean inversion, if the sphere is a Euclidean sphere, or a Euclidean reflection, if the sphere is a Euclidean plane. The group of \emph{M\"obius transformations} on $\mathbb{R}^N_\infty$ is the group generated by inversions in spheres.  By an \emph{open ball} $B$ in $\mathbb{R}^N_{\infty}$ we mean one of the connected components of the complement of a sphere. We write $\widehat{B}$ for the other component. Given two distinct open balls $B_1$ and $B_2$ in $\mathbb{R}^N_\infty$, we denote the signed inversive distance between them by $[B_1, B_2]$. For any four distinct points $p_1$, $p_2$, $p_3$, and $p_4$ in $\mathbb{R}^N_\infty$ we let $|p_1, p_2, p_3, p_4|$ denote their absolute cross-ratio. The signed inversive distance and absolute cross-ratio are two geometric quantities invariant under the action of the M\"{o}bius group on $\mathbb{R}^N_\infty$. We will describe them in detail in section~\ref{S: background}. The two main results of this paper state that these invariants suffice to rigidify a configuration of open balls, or a configuration of points, up to M\"{o}bius transformations. 

\begin{theorem}\label{T: main1}
Let $\{B_\alpha : \alpha \in \mathcal{A}\}$ and $\{B'_\alpha : \alpha \in \mathcal{A}\}$ be two collections of open balls in $\RN$, indexed by the same set. Suppose that $\bigcap_{\alpha\in \mathcal{A}} \partial B_\alpha=\emptyset$. Then there is a M\"obius transformation $f$ such that one of the following holds: either $f(B_\alpha)=B'_\alpha$ for each $\alpha$ in $\mathcal{A}$, or else $f(\widehat{B_\alpha})=B'_\alpha$ for each $\alpha$ in $\mathcal{A}$, if and only if $[B_\alpha,B_\beta]=[B'_\alpha,B'_\beta]$ for all pairs  $\alpha$ and $\beta$ in $\mathcal{A}$.
\end{theorem}





\begin{theorem}\label{T: main2}
Let $\{p_\alpha\,:\,\alpha\in \mathcal{A}\}$ and $\{p'_\alpha\,:\,\alpha\in \mathcal{A}\}$ be two collections of distinct points in $\RN$, indexed by the same set.  There is a M\"obius transformation $f$ with $f(p_\alpha)=p'_\alpha$ for each $\alpha$ in $\mathcal{A}$ if and only if $|p_\alpha,p_\beta,p_\gamma,p_\delta|=|p'_\alpha,p'_\beta,p'_\gamma,p'_\delta|$ for all ordered $4$-tuples $(\alpha, \beta, \gamma, \delta)$ of distinct indices in $\mathcal{A}$.
\end{theorem}

These theorems resolve two problems posed by Beardon and Minda \cite{BeMi2006} concerning extensions and 
higher-dimensional generalizations of their results on the conformal symmetries of circular regions in the extended complex plane.

\section{Background}\label{S: background}

 A \emph{circular region} in the extended complex plane  $\mathbb{C}_\infty$ is a region bounded by a collection of pairwise disjoint circles. A \emph{finitely connected region} in $\mathbb{C}_\infty$ is a region with a finite number of boundary components. A classical theorem of Koebe (that can be found in \cite[chapter 15]{Co1995} or \cite[chapter X]{Fo1951}) says that a finitely connected region is conformally equivalent to a finitely connected circular region that has a finite number of punctures. A \emph{M\"obius transformation} is a conformal or anti-conformal homeomorphism of $\mathbb{C}_\infty$. Such maps can be expressed algebraically as
\begin{equation}\label{E: Mobius}
z\mapsto \frac{az+b}{cz+d},\quad \text{or} \quad z\mapsto \frac{a\bar{z}+b}{c\bar{z}+d},
\end{equation}
where $ad-bc\neq 0$.  Given two Euclidean circles $C_1$ and $C_2$, with centres $c_1$ and $c_2$, and radii $r_1$ and $r_2$, the \emph{inversive distance} between these two circles is the positive quantity
\begin{equation}\label{E: inversive1}
(C_1,C_2) = \left|\frac{r_1^2+r_2^2-|c_1-c_2|^2}{2r_1r_2}\right|.
\end{equation} 
More generally, if $C_1$ and $C_2$ are two circles in $\mathbb{C}_\infty$ (that is, they are each either Euclidean circles, or Euclidean lines with the point $\infty$ attached), then we define the inversive distance $(C_1,C_2)$ to be $(f(C_1),f(C_2))$, where $f$ is any M\"obius transformation that maps both $C_1$ and $C_2$ to Euclidean circles. This definition is independent of $f$, and the resulting quantity is invariant under M\"obius transformations, in the sense that $(g(C_1),g(C_2))=(C_1,C_2)$ for each M\"obius map $g$. See \cite[section 3.2]{Be1983} for information on the inversive distance.

The following result is part of \cite[Thm 4.1]{BeMi2006}; the original theorem of Beardon and Minda also includes a uniqueness statement, which we will address in section~\ref{S: proof1}. 

\newtheorem*{thmA}{Theorem A}
\begin{thmA}
Suppose that $\Omega$ and $\Omega'$ are circular regions bounded by circles $C_1,\dots,C_m$ and $C'_1,\dots,C'_m$, respectively, where $m\geq 2$. Then there is a M\"obius transformation $f$ with $f(\Omega)=\Omega'$ and $f(C_j)=C'_j$, $1\leq j\leq m$, if and only if $(C_j,C_k)=(C'_j,C'_k)$ for all $j$ and $k$ with $1\leq j < k \leq m$.
\end{thmA}

Beardon and Minda also gave an analogous result about punctured regions.  For points $a$, $b$, $c$, and $d$ in $\mathbb{C}_\infty$, let $|a,b,c,d|$ denote the \emph{absolute cross-ratio} of $a$, $b$, $c$, and $d$; that is, 
\begin{equation}\label{E: cross-ratio}
|a,b,c,d| = \frac{|a-b||c-d|}{|a-c||b-d|},
\end{equation}
with the usual conventions regarding the point $\infty$.
The absolute cross-ratio is invariant under M\"obius transformations. The following result is \cite[Theorem 14.1]{BeMi2006}.

\newtheorem*{thmB}{Theorem B}
\begin{thmB}
Given two collections of points $p_1,\dots,p_m$ and $p'_1,\dots,p'_m$ in $\mathbb{C}_\infty$, $m\geq 4$, there is a M\"obius transformation $f$ with $f(p_i)=p'_i$ for $i=1,2,\dots,m$ if and only if $|p_i,p_j,p_k,p_l|=|p'_i,p'_j,p'_k,p'_l|$ for all $1\leq i,j,k,l\leq m$.
\end{thmB}

A weaker theorem than Theorem B in which the absolute cross-ratio is replaced by the usual complex cross-ratio is well known and straightforward to prove.

At the end of \cite{BeMi2006}, Beardon and Minda asked the following questions (the second question has been paraphrased).

\theoremstyle{definition}
\newtheorem*{Q1}{Question 1}
\begin{Q1}
Is the conclusion of Theorem A valid when $C_1,\dots,C_m$ are any set of $m$ distinct circles in $\mathbb{C}_\infty$? Here the $C_i$ are allowed to be intersecting, or tangent, to each other.
\end{Q1}

\newtheorem*{Q2}{Question 2}
\begin{Q2}
Do Theorems A and B generalize to higher dimensions?
\end{Q2}

The answer to Question 1 is negative, and we provide examples to justify this in section \ref{S: examples}. Subject to certain restrictions, however, both Theorems A and B generalize to allow arbitrarily many circles and points, and the circles may intersect. These generalizations are our main theorems (Theorems~\ref{T: main1} and \ref{T: main2}), and they apply in all dimensions.

To generalize Theorem A we work with the \emph{signed} inversive distance between \emph{discs} rather than circles (or, in higher dimensions, with balls rather than spheres). Given two Euclidean balls $B_1$ and $B_2$, with centres $c_1$ and $c_2$, and radii $r_1$ and $r_2$, the \emph{signed inversive distance} between these two balls is the quantity
\begin{equation}\label{E: inversive2}
[B_1,B_2] = \frac{r_1^2+r_2^2-|c_1-c_2|^2}{2r_1r_2}.
\end{equation} 
Again, $[B_1,B_2]$ can be defined for arbitrary balls $B_1$ and $B_2$ by transferring away from the point $\infty$ using a M\"obius transformation, and again, the quantity $[B_1,B_2]$ is preserved under M\"obius transformations. Notice that $[\widehat{B_1},B_2]=-[B_1,B_2]$ and $(\partial B_1,\partial B_2)=|[B_1,B_2]|$. 



To recover Theorem A from the $N=2$ case of Theorem~\ref{T: main1}, we begin with the hypotheses of Theorem A, and define $B_i$ to be the component of $\mathbb{C}_\infty\setminus C_i$ that contains $\Omega$. Likewise we define $B'_i$ to be the component of $\mathbb{C}_\infty\setminus C'_i$ that contains $\Omega'$. This means that $[B_i,B_j]=(C_i,C_j)$ and $[B'_i,B'_j]=(C'_i,C'_j)$ for all indices $i$ and $j$ in $\{1,\dots,m\}$. From Theorem~\ref{T: main1} we deduce the existence of a M\"obius map $f$ such that one of the following holds: either $f(B_i)=B_i'$ for each $i$, or else $f(\widehat{B_i})=B_i'$ for each $i$. In the latter case, because $\widehat{B_2}\subset B_1$, we find that 
\[
B_2'=f(\widehat{B_2})\subset f(B_1)=\widehat{B_1'},
\]
which is false. Therefore $f(B_i)=B_i'$ for each $i$, which means that $f(C_i)=C_i'$ for each $i$, and because $\Omega =\bigcap_{i=1}^m B_i$ and $\Omega' =\bigcap_{i=1}^m B'_i$, we also see that $f(\Omega)=\Omega'$.

Theorem~\ref{T: main1} may fail  when $\bigcap_{\alpha\in \mathcal{A}} \partial B_\alpha=\emptyset$; an example of its failure is given in section \ref{S: examples}. 

\section{Examples}\label{S: examples}

We provide a sequence of examples which answer Question 1 and explain the necessity of the conditions in Theorem~\ref{T: main1}.

\begin{figure}[ht]
\centering
\includegraphics[scale=1.0]{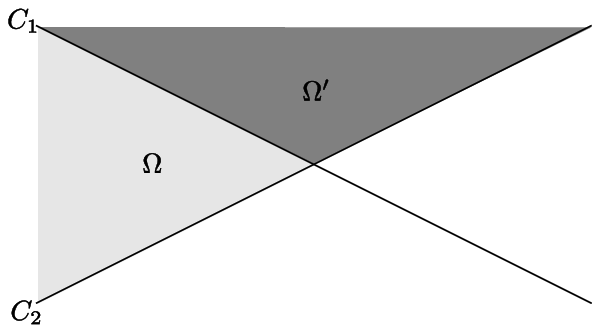}
\caption{}
\label{F: example4}
\end{figure}

\begin{example}\label{E: easy}
Here is the simplest example to show that Theorem A is invalid when the circles $C_i$ are allowed to intersect.
Let $C_1$ and $C_2$ be two Euclidean lines through the origin that cross at an angle $\pi/3$. Let $C_1'=C_1$ and $C_2'=C_2$. Define $\Omega$ to be one of the resulting sectors with angle $\pi/3$, and define $\Omega'$ to be one of the sectors with angle $2\pi/3$. See  Figure~\ref{F: example4}. Only M\"obius transformations of the form  $z\mapsto \lambda z$ or $z\mapsto \lambda/z$, for non-zero real numbers $\lambda$, fix both $C_1$ and $C_2$ as sets. None of these transformations map $\Omega$ to $\Omega'$. 
\end{example}

There is not a M\"obius transformation mapping the set $\Omega$ in Example~\ref{E: easy} to $\Omega'$; however, there is a M\"obius transformation that maps $C_1$ to $C'_1$ and $C_2$ to $C'_2$, namely the identity. This explains why we consider balls rather than spheres in Theorem~\ref{T: main1}.

\begin{figure}[ht]
\centering
\includegraphics[scale=1.0]{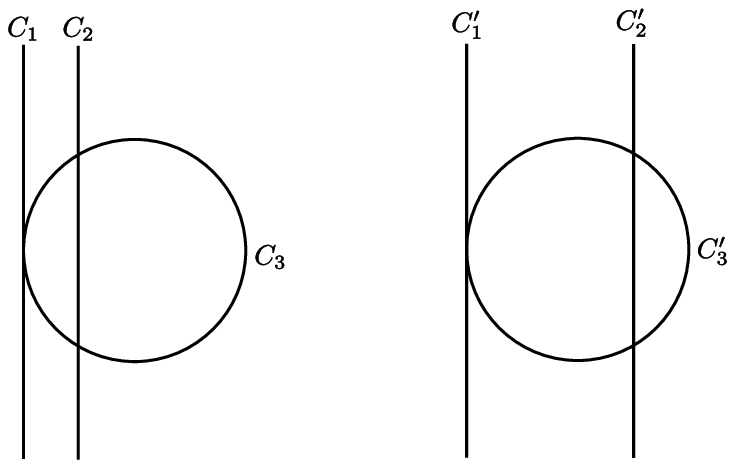}
\caption{}
\label{F: example5}
\end{figure}

\begin{example}\label{E: simple}
This example shows that it is necessary to use signed inversive distances. Let $C_1$ and $C_1'$ both denote the line $x=-1$. Let $C_2$ denote the line $x=-1/2$ and let $C'_2$ denote the line $x=1/2$. Let $C_3$ and $C'_3$ both denote the unit circle. See Figure~\ref{F: example5}. Then $(C_i,C_j)=(C'_i,C'_j)$ for all $1\leq i,j\leq 3$, but there is not a M\"obius transformation that maps $C_i$ to $C_i'$ for $i=1,2,3$.
\end{example}

The inversive distance condition is not sensitive enough to distinguish the two geometric configurations shown in Figure~\ref{F: example5}. It is impossible to define discs $B_i$ and  $B_i'$ with boundary circles $C_i$ and $C_i'$ in Example~\ref{E: simple} such that $[B_i,B_j]=[B'_i,B'_j]$ for all $i$ and $j$ in $\{1,2,3\}$. 

The final example in this section shows that Theorem~\ref{T: main1} would fail if we allowed $C_1\cap\dots\cap C_m \neq \emptyset$.

\begin{figure}[ht]
\centering
\includegraphics[scale=1.0]{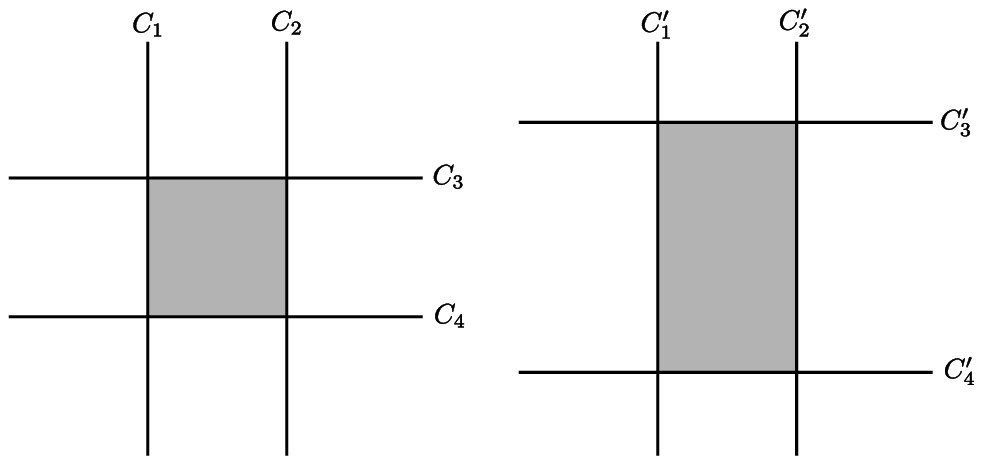}
\caption{}
\label{F: rectangle}
\end{figure}

\begin{example}\label{E: rectangle}
Let $C_1$, $C_2$, $C_3$, and $C_4$ be the extended sides of a square in the complex plane. Let $C'_1$, $C'_2$, $C'_3$, and $C'_4$ be the extended sides of a rectangle (that is not a square) in the complex plane. See Figure~\ref{F: rectangle}. For each $i$, let $B_i$ be the half-plane with boundary $C_i$ that contains the shaded square, and let $B_i'$ be the half-plane with boundary $C'_i$ that contains the shaded rectangle. Then $[B_i,B_j]=[B'_i,B'_j]$ for all pairs $i,j$. If there is a M\"obius transformation $f$ that maps $B_i$ to $B_i'$ for each $i$, then  $f$ must map the square to the rectangle. This cannot be. 
\end{example}

\section{Hyperbolic geometry}\label{S: hyperbolic}

Beardon and Minda noted that their results can be interpreted in terms of hyperbolic geometry,  and this is our starting point. Refer to  \cite{Be1983,Ra1994} for complete introductions to hyperbolic geometry. 

The action of M\"obius transformations on $\mathbb{C}_\infty$ extends to an action on the upper half-space model of three-dimensional  hyperbolic space, $\mathbb{H}^3$, and this action on $\mathbb{H}^3$ is isometric with respect to the hyperbolic metric, $\rho$. Each circle in $\mathbb{C}_\infty$ is the ideal boundary of a unique hyperbolic plane in $\mathbb{H}^3$. If $\Pi_1$ and $\Pi_2$ are two hyperbolic planes with ideal boundary circles $C_1$ and $C_2$, then
\begin{equation*}\label{E: inversive3}
(C_1,C_2) = 
\begin{cases}
\cosh \rho(\Pi_1,\Pi_2) & \text{if $\Pi_1$ and $\Pi_2$ are disjoint},\\
\cos \theta & \text{if $\Pi_1$ and $\Pi_2$ intersect in an angle $\theta$}.
\end{cases}
\end{equation*}
In higher dimensions, the situation is similar. The set
\[
\left\{ (x_1,\dots,x_{N+1})\in\mathbb{R}^{N+1}\,:\, x_{N+1}>0\right\},
\]
which we denote by $\mathbb{H}^{N+1}$, is a model of $(N+1)$-dimensional hyperbolic space when equipped with the Riemannian density $ds = |dx|/x_{N+1}$. The ideal boundary of $\mathbb{H}^{N+1}$ consists of the set $x_{N+1}=0$, which we identify with $\mathbb{R}^N$, and the point $\infty$. The action of the group of M\"obius transformations on $\mathbb{R}^N_\infty$ extends to an isometric action on $\mathbb{H}^{N+1}$. The boundary in $\mathbb{H}^{N+1}$ of a hyperbolic half-space $\Sigma$ is an $N$-dimensional hyperbolic plane $\Pi$.  The ideal boundary of $\Sigma$ is a spherical ball $B$ in $\mathbb{R}^{N}_\infty$. The boundary in $\mathbb{R}^{N}_\infty$ of $B$ is a sphere $C$, and $C$ is the ideal boundary of $\Pi$. We often move between spheres and balls in $\mathbb{R}^{N}_\infty$, and half-spaces and planes in $\mathbb{H}^{N+1}$. 

There are many models of $(N+1)$-dimensional hyperbolic space, and although $\mathbb{H}^{N+1}$ is the most appropriate model for explaining how Theorem~\ref{T: main1} can be considered as a problem in hyperbolic geometry, the \emph{hyperboloid model} of hyperbolic space is the most appropriate model for proving the theorem. We describe the hyperboloid model in the next section. 

\section{The hyperboloid model of hyperbolic space}\label{S: hyperboloid}

The substance of this section is taken from \cite[chapter 3]{Ra1994}.

We equip $\mathbb{R}^{N+1}$ with the \emph{Lorentz inner product} $\langle\,\cdot\,\rangle$, defined by 
\[
\la(x_1,\dots,x_{N+1}),(y_1,\dots,y_{N+1})\ra = x_1y_1+\dots + x_Ny_N-x_{N+1}y_{N+1}.
\]
This is not an inner product in the usual sense, as it is not positive definite. We write $\|x\|^2 = \la x,x\ra$.  In contrast, we denote the Euclidean scalar product of points $x$ and $y$ in $\mathbb{R}^{N+1}$ by $x\cdot y$, and the Euclidean norm of $x$ by $|x|$. A vector $x$ in $\mathbb{R}^{N+1}$ is \emph{space-like} if $\|x\|^2>0$, \emph{time-like} if $\|x\|^2<0$, and \emph{light-like} if $\|x\|^2=0$. The terminology originates from the theory of relativity. A subspace $V$ of $\mathbb{R}^{N+1}$ is \emph{space-like} if every non-zero element of $V$ is space-like, \emph{time-like} if there is a time-like vector in $V$, and \emph{light-like} otherwise.

A linear map of $\mathbb{R}^{N+1}$ that preserves the Lorentz inner inner product is described as a \emph{Lorentz transformation}. A vector $x$ in $\mathbb{R}^{N+1}$ is \emph{positive} if $x_{N+1}>0$. A Lorentz transformation is \emph{positive} if it maps positive time-like vectors to positive time-like vectors. The positive Lorentz transformations together form a group, denoted $\textup{PO}(N,1)$.

The underlying space of the hyperboloid model of $N$-dimensional hyperbolic space is the hyperboloid sheet 
\[
\mathcal{H}^N = \{\, x\in\mathbb{R}^{N+1}\,:\|x\|^2=-1,\,x_{N+1}>0\,\},
\]
embedded in $\mathbb{R}^{N+1}$. This is a model of $N$-dimensional hyperbolic space with the metric $\rho$ defined by
\[
\cosh \rho(x,y) =-\langle x,y\rangle.
\] 

The group $\textup{PO}(N,1)$ consists of those Lorentz transformations that fix $\mathcal{H}^N$ (the Lorentz transformations that are not positive swap $\mathcal{H}^N$ with its twin hyperboloid sheet). The group $\textup{PO}(N,1)$ is the full group of hyperbolic isometries of $\mathcal{H}^N$.

The hyperbolic lines in $\mathcal{H}^N$ are intersections of $\mathcal{H}^N$ with two-dimensional time-like subspaces of $\mathbb{R}^{N+1}$. The hyperbolic planes of codimension one in $\mathcal{H}^N$ are intersections of $\mathcal{H}^N$ with $N$-dimensional time-like subspaces of $\mathbb{R}^{N+1}$. In future we describe such planes merely as `planes', because all the planes we consider have codimension one.


Given a subspace $V$ of $\mathbb{R}^{N+1}$, the \emph{Lorentz complement} of $V$ is the space
\[
V^{\textup{L}} = \{\,y\in\mathbb{R}^{N+1}\,:\, \langle x,y\rangle=0\text{ for all $x$ in $V$}\,\}.
\]
 To each time-like Euclidean plane $P$ there corresponds a unique line $\ell$ of space-like vectors in $\mathbb{R}^{N+1}$ that are Lorentz orthogonal to $P$, so that $P=\ell^{\textup{L}}$. Conversely, to a Euclidean line $\ell$ of space-like vectors there corresponds a unique time-like Euclidean plane $P$ that is Lorentz orthogonal to $\ell$.

Let $P_1$ and $P_2$ be two $N$-dimensional time-like planes in $\mathbb{R}^{N+1}$ with non-zero space-like normals $v_1$ and $v_2$, respectively, where $\|v_1\|^2=\|v_2\|^2=1$. Let $\Pi_1=\mathcal{H}^N\cap P_1$ and $\Pi_2=\mathcal{H}^N\cap P_2$. The ideal boundaries of $\Pi_1$ and $\Pi_2$ are spheres $C_1$ and $C_2$. The inversive distance of $C_1$ and $C_2$ defined in \eqref{E: inversive1} satisfies
\begin{equation}\label{E: inversiveProd2}
(C_1,C_2) = |\langle v_1,v_2\rangle |
\end{equation}
(see \cite[section 3.2]{Ra1994}). This is the simplest formula for the inversive distance so far, hinting that the hyperboloid model may be the most natural setting for considering Theorem~\ref{T: main1}. 

The plane $P_1$ consists of all points $x$ in $\mathbb{R}^{N+1}$ for which $\langle x,v_1\rangle=0$. Let $Q_1$ consist of all points $x$ in $\mathbb{R}^{N+1}$ for which $\langle x,v_1\rangle>0$. Define the half-space $\Sigma_1$ to be equal to $\mathcal{H}^N\cap Q_1$. We define $\Sigma_2$ in a similar fashion using $P_2$. The ideal boundaries of $\Sigma_1$ and $\Sigma_2$ are open spherical balls $B_1$ and $B_2$. The signed inversive distance of $B_1$ and $B_2$ defined in \eqref{E: inversive2} satisfies the formula
\begin{equation}\label{E: inversiveProd}
[B_1,B_2] = \langle v_1,v_2\rangle. 
\end{equation}
(again, see \cite[section 3.2]{Ra1994}).

\section{Canonical forms for subspaces of Lorentz space}\label{S: canonical}

Let $e_1,\dots,e_{N+1}$ be the standard basis vectors. For each $p=1,\dots,N$, define subspaces 
\begin{align*}
T_p &=\{\,(x_1,\dots,x_{p-1},0,\dots,0,x_{N+1})\in\mathbb{R}^{N+1}\,:\,x_i\in\mathbb{R}\,\}, \\
S_p &=\{\,(x_1,\dots,x_{p-1},x_{p},0,\dots,0,0)\in\mathbb{R}^{N+1}\,:\,x_i\in\mathbb{R}\,\}, \\
L_p &=\{\,(x_1,\dots,x_{p-1},\lambda,0,\dots,0,\lambda)\in\mathbb{R}^{N+1}\,:\,\lambda,x_i\in\mathbb{R}\,\},
\end{align*}
each of dimension $p$. We identify the subspace $S_p$  with $\mathbb{R}^p$, for each $p$. Notice that $T_p$ is time-like, since it contains the time-like vector $e_{N+1}$; $S_p$ is space-like, since each non-zero vector in $S_p$ is space-like; and $L_p$ is light-like, because it contains no time-like vectors, but it does contain the light-like vector $e_p+e_{N+1}$. 
\begin{lemma}\label{L: canonical}
Each $p$-dimensional subspace of $\mathbb{R}^{N+1}$ (for $1 \leq p \leq N$) is isomorphic by a Lorentz transformation to either $T_p$, $S_p$, or $L_p$.
\end{lemma}
\begin{proof}
Given a $p$-dimensional proper subspace $V$, let $\alpha$ be a Lorentz transformation that fixes $e_{N+1}$ and acts as a standard orthogonal map on $\mathbb{R}^N$ in such a way that $\mathbb{R}^N\cap V$ is mapped to $\mathbb{R}^k$, where $k$ is the dimension of $\mathbb{R}^N\cap V$. Either $k=p$,  in which case $V$ is contained in $\mathbb{R}^N$ and the proof is finished, or $k=p-1$. In the second case, choose an element $u$ in $\alpha(V)\setminus\mathbb{R}^{N}$. Let 
\[
v = (0,\dots,0,u_p,\dots,u_{N+1});
\]
this vector is also in $\alpha(V)\setminus\mathbb{R}^{N}$, since $\mathbb{R}^{p-1}$ is contained in $\alpha(V)$. Choose a Lorentz transformation $\beta$ that fixes $e_1,\dots,e_{p-1}$ and $e_{N+1}$, and acts as a standard orthogonal map on the span of $e_p,\dots,e_N$ in such a way that $v$ maps to
\[
w = (0,\dots,0,w_p,0\dots,0,w_{N+1}),
\]
where $w_{N+1}=u_{N+1}$. Note that $\beta\alpha(V)$ is the span of $\mathbb{R}^{p-1}$ and $w$. Let $A[a,b]$ denote the $(N+1)$-by-$(N+1)$ Lorentz matrix whose entries $A_{i,j}$ coincide with the entries of the identity matrix, except $A_{p,p}=A_{N+1,N+1}=a$ and $A_{p,N+1}=A_{N+1,p}=b$, where $a^2-b^2=1$. We define a third Lorentz transformation $\gamma$ as follows. If $w_p=w_{N+1}$ then let $\gamma$ be the identity map, and if $w_p=-w_{N+1}$ then let $\gamma$ be the map $(x_1,\dots,x_N,x_{N+1})\mapsto (x_1,\dots,x_N,-x_{N+1})$. Otherwise, let $\ell(w) = \sqrt{\left|w_p^2-w_{N+1}^2\right|}$ and define
\[
\gamma = 
\begin{cases}
A\left[\tfrac{w_p}{\ell(w)},-\tfrac{w_{N+1}}{\ell(w)}\right] & \text{if $\|w\|^2>0$,}\vspace{1mm}\\
A\left[\tfrac{w_{N+1}}{\ell(w)},-\tfrac{w_{p}}{\ell(w)}\right] & \text{if $\|w\|^2<0$.}
\end{cases}
\]
 The Lorentz transformation $\gamma\beta\alpha$ maps $V$ to  either $L_p$, $S_p$, or $T_p$, depending on whether $\|w\|^2=0$, $\|w\|^2>0$, or $\|w\|^2<0$.
\end{proof}

\section{Proof of Theorem~\ref{T: main1}}\label{S: proof1}

The proofs of both Theorem~\ref{T: main1} and Theorem~\ref{T: main2} are based on the following proposition.

\begin{proposition}\label{P: existence}
Let $\{v_\alpha\,\,:\,\alpha\in \mathcal{A}\}$ and $\{v'_\alpha\,\,:\,\alpha\in \mathcal{A}\}$ be two collections of vectors in $\mathbb{R}^{N+1}$ such that $\langle v_\alpha,v_\beta\rangle = \langle v'_\alpha,v'_\beta\rangle$ for all pairs $\alpha$ and $\beta$ in $\mathcal{A}$. Suppose that the subspace  spanned by the $v_\alpha$ is either time-like or space-like. Then there is a  Lorentz transformation $\phi$ with $\phi(v_\alpha)=v'_\alpha$ for each $\alpha$ in $\mathcal{A}$.
\end{proposition}

Proposition~\ref{P: existence} fails when the subspace $V$ spanned by the $v_\alpha$ is light-like because  the next elementary lemma, used in the proof of Proposition~\ref{P: existence}, also fails when $V$ is light-like. 

\begin{lemma}\label{L: failure}
Let $V$ be either a time-like or a space-like subspace of $\mathbb{R}^{N+1}$. If there is an element $v$ of $V$ such that $\langle v,w\rangle = 0$ for all vectors $w$ in $V$ then $v=0$.
\end{lemma}

Indeed, if $\langle v,w\rangle = 0$ for all $w$ in $V$ then, in particular,   $\langle v, v \rangle = 0$ so $v$ is either $0$ or light-like. 

\begin{proof}[Proof of Proposition~\ref{P: existence}]
Let $V$ denote the subspace spanned by the vectors $v_\alpha$, and let $V'$ denote the subspace spanned by the vectors $v'_\alpha$.
By applying preliminary Lorentz transformations, we may assume that each of $V$ and $V'$ are either equal to $\mathbb{R}^{N+1}$ or else assume one of the canonical forms listed at the beginning of section \ref{S: canonical}. The subspace $V'$ is time-like if $V$ is time-like, and $V'$ is space-like if $V$ is space-like. The following proof is valid whether $V$ is time-like or space-like.

Let $p$ and $q$ be the dimensions of $V$ and $V'$. By swapping $V$ and $V'$ if necessary we may assume that $p\geq q$. Relabel the $v_\alpha$ so that  $v_1,\dots,v_p$ span $V$, and the remaining vectors $v_\alpha$ are linearly dependent on $v_1,\dots,v_p$. Relabel the $v'_\alpha$ in a corresponding fashion.
It will now be shown that $v'_1,\dots,v'_p$ are linearly independent. Suppose that $\lambda_1 v'_1+\dots +\lambda_p v'_p=0$, for real numbers $\lambda_1,\dots,\lambda_p$. Define $v=\lambda_1 v_1+\dots +\lambda_p v_p$. Then $\langle v,v_i\rangle= \langle 0,v'_i\rangle=0$, for $i=1,\dots,p$. Thus $v=0$, by Lemma~\ref{L: failure}. By linear independence of $v_1,\dots,v_p$, we deduce that the numbers $\lambda_i$ are all $0$. We can extend both $v_1,\dots,v_p$ and $v'_1,\dots,v'_p$ to bases of $\mathbb{R}^{N+1}$ using either the standard basis vectors $e_{p+1},\dots,e_{N+1}$ (if $V=S_p$) or $e_p,\dots,e_N$ (if $V=T_p$).

Let $\phi$ be the unique bijective linear map that fixes each of these $N+1-p$ standard basis vectors, and satisfies $\phi(v_j)=v'_j$, for $j=1,\dots,p$. Because the vectors $v_i$ and $v'_j$ are Lorentz orthogonal to the vectors $e_k$ the map $\phi$ is a Lorentz transformation. 

Finally, observe that for indices $\alpha$ other than $1,\dots,p$, we have \[\langle\phi(v_\alpha)-v'_\alpha,v'_j\rangle=\langle v_\alpha,v_j\rangle-\langle v'_\alpha,v'_j\rangle=0, \text{\;for\;} j=1,\dots,p.\] Therefore $\phi(v_\alpha)=v'_\alpha$, by Lemma~\ref{L: failure}. 
\end{proof}

Before we prove Theorem~\ref{T: main1}, we state a lemma that explains the significance of the condition $\bigcap_{\alpha\in\mathcal{A}}\partial B_\alpha = \emptyset$ of Theorem~\ref{T: main1} in terms of hyperbolic geometry. (Note that, because the ideal boundary of $\mathcal{H}^N$ is $(N-1)$-dimensional, we assume that $B_\alpha$ and $B'_\alpha$ are balls in $\mathbb{R}^{N-1}_\infty$, rather than $\mathbb{R}^N_\infty$.) The spheres $\partial B_\alpha$ are the ideal boundaries of hyperbolic planes $\Pi_\alpha$, and each hyperbolic plane $\Pi_\alpha$ is the intersection of $\mathcal{H}^N$ with a time-like Euclidean plane $P_\alpha$ in $\mathbb{R}^{N+1}$. Let $\Sigma_\alpha$ denote the hyperbolic half-space with ideal boundary $B_\alpha$, and let  $v_\alpha$ denote the unique space-like Lorentz unit normal of $P_\alpha$ such that $\Sigma_\alpha = \{x\,:\, \langle x,v_\alpha\rangle >0\}$.

\begin{lemma}\label{L: noIntersection}
The spheres $\partial B_\alpha$ do not contain a common point of intersection if and only if the subspace of $\mathbb{R}^{N+1}$ spanned by the $v_\alpha$ is either time-like, or $N$-dimensional and space-like.
\end{lemma}
\begin{proof}
Let $V$ be the subspace spanned by the vectors $v_\alpha$. By applying a Lorentz transformation we can assume that $V$ is either $\mathbb{R}^{N+1}$, $T_p$, $S_p$, or $L_p$. The spheres $\partial B_\alpha$ contain a common point of intersection if and only if there is a light-like vector in $\bigcap_{\alpha\in \mathcal{A}} P_\alpha$.  That is, if and only if there is a light-like vector in $V^{\textup{L}}$. If $V=L_p$ then the vector $e_p-e_{N+1}$ is light-like and contained in $V^{\textup{L}}$, and if $V=S_p$ and $p<N$ then the vector $e_N+e_{N+1}$ is light-like and contained in $V^{\textup{L}}$. If $V=\mathbb{R}^{N+1}$, $V=T_p$, or $V=S_N$, then $V^{\textup{L}}$ does not contain any light-like vectors.
\end{proof}

 Let $\widehat{\Sigma_\alpha}$ denote the hyperbolic half-space with ideal boundary $\widehat{B_\alpha}$.

\begin{proof}[Proof of Theorem~\ref{T: main1}]
Suppose that $[B_\alpha,B_\beta]=[B'_\alpha,B'_\beta]$ for each $\alpha$ and $\beta$ in $\mathcal{A}$. We must construct a M\"obius transformation $f$ for which either $f(B_\alpha)=B'_\alpha$ for all $\alpha$, or else $f(\widehat{B_\alpha})=B'_\alpha$ for all $\alpha$. (Note that the converse implication follows immediately by preservation of the signed inversive distance under M\"obius transformations.)

Since $\bigcap_{\alpha\in \mathcal{A}}\partial B_\alpha=\emptyset$, we see from Lemma~\ref{L: noIntersection} that the subspace $V$ spanned by the $v_\alpha$ is not light-like. Proposition~\ref{P: existence} shows that there is a Lorentz transformation $\phi$ such that $\phi(\Sigma_\alpha)=\Sigma'_\alpha$ for each $\alpha$ in $\mathcal{A}$. Either $\phi$ or $-\phi$ is a positive Lorentz transformation. In the first case $\phi(\Sigma_\alpha)=\Sigma'_\alpha$ for each $\alpha$ in $\mathcal{A}$, and in the second case $-\phi(\widehat{\Sigma_\alpha})=\Sigma'_\alpha$ for each $\alpha$ in $\mathcal{A}$. The action of this positive Lorentz transformation on the ideal boundary of hyperbolic space is a M\"obius transformation with the required properties.
\end{proof}

The map $f$ of Theorem~\ref{T: main1} is uniquely determined if and only if the map $\phi$ of Proposition~\ref{P: existence} is uniquely determined. This occurs if and only if the subspace $V$ spanned by the $v_\alpha$ is the whole of $\mathbb{R}^{N+1}$. 

From Lemma~\ref{L: noIntersection} we know that $V\neq \mathbb{R}^{N+1}$ if and only if either $V$ is (i) $N$-dimensional and space-like, or (ii) time-like of dimension less than $N+1$. Case (i) occurs if and only if $V^\textup{L}$ is $1$-dimensional and time-like. Since $\bigcap_{\alpha\in\mathcal{A}}P_\alpha = V^{\textup{L}}$ and $\bigcap_{\alpha\in\mathcal{A}}\Pi_\alpha = \mathcal{H}^N\cap \left(\bigcap_{\alpha\in\mathcal{A}}P_\alpha\right)$, case (i) occurs if and only if the intersection of the planes $\Pi_\alpha$ is a single point in hyperbolic space. Case (ii) occurs if and only if there is a space-like vector that is Lorentz orthogonal to $V$. That is, if and only if there is a time-like Euclidean plane that is Lorentz orthogonal to all the $P_\alpha$. Or, equivalently, there is a sphere in $\mathbb{R}^{N-1}_\infty$ that is orthogonal to all the $\partial B_\alpha$.

Since Beardon and Minda considered only non-intersecting circles, case (i) did not arise in their study. Case (ii) did arise: they defined a collection of circles in $\mathbb{C}_\infty$ to be \emph{strongly symmetric} if there is another circle orthogonal to each circle in the collection. Beardon and Minda prove, as we have just verified, that the map $f$ is unique if and only if the collection of circles is not strongly symmetric.

\section{Proof of Theorem~\ref{T: main2}}\label{S: proof2}

Proposition~\ref{P: existence}  can also be used to  prove Theorem~\ref{T: main2}. To apply this proposition, first the $(N-1)$-dimensional unit sphere must be identified with the ideal boundary of $\mathcal{H}^N$. The $N$-dimensional unit ball $\mathbb{B}^{N}$ is also a model of hyperbolic space, and there is an isometry $\Phi$ from $\mathbb{B}^N$ to $\mathcal{H}^{N}$ given by 
\[
(x_1,\dots,x_N)\mapsto \left(\frac{2x_1}{1-|x|^2},\dots,\frac{2x_N}{1-|x|^2},\frac{1+|x|^2}{1-|x|^2}\right).
\]
With this correspondence, the point $(x_1,\dots,x_N)$ in $\mathbb{S}^{N-1}$ is paired with the Euclidean line that passes through $0$ and $(x_1,\dots,x_N,1)$. The absolute cross-ratio of four points $p_1$, $p_2$, $p_3$, and $p_4$ in $\mathbb{S}^{N-1}$ is 
\begin{equation}\label{E: ratio1}
|p_1,p_2,p_3,p_4|=\frac{|p_1-p_2||p_3-p_4|}{|p_1-p_3||p_2-p_4|}.
\end{equation}
We now wish to define the cross-ratio in terms of the Lorentz model of hyperbolic space. We use the next elementary lemma.

\begin{lemma}
If $u$ and $v$ are two linearly independent  positive light-like vectors in $\mathbb{R}^{N+1}$ then $\langle u,v\rangle<0$. 
\end{lemma}
\begin{proof}
Since $u$ and $v$ are positive and light-like we can choose elements $u_0$ and $v_0$ of $\mathbb{R}^N$ such that $u=u_0+|u_0|e_{N+1}$ and $v=v_0+|v_0|e_{N+1}$. Therefore $|u|=\sqrt{2}|u_0|$ and $|v|=\sqrt{2}|v_0|$. We obtain
\[
\langle u,v\rangle = u_0\cdot v_0 - |u_0||v_0| = u\cdot v - 2|u_0||v_0| < |u||v|-2|u_0||v_0|=0
\]
by the Cauchy-Schwarz inequality.
\end{proof}

Given light-like lines $\ell_1$, $\ell_2$, $\ell_3$, and $\ell_4$ we choose, for each $i$, any positive light-like vector $v_i$ in $\ell_i$ and define
\begin{equation*}\label{E: ratio2}
|\ell_1,\ell_2,\ell_3,\ell_4|=\frac{\langle v_1,v_2\rangle \langle v_3,v_4\rangle}{\langle v_1,v_3\rangle \langle v_2,v_4\rangle}.
\end{equation*}
This quantity is preserved under Lorentz transformations. If the point $p_i$ in $\mathbb{S}^{N-1}$ corresponds to the line $\ell_i$ under $\Phi$ then we have
\[
|\ell_1,\ell_2,\ell_3,\ell_4|=|p_1,p_2,p_3,p_4|^2.
\]
It suffices to verify this formula when $p_1=e_1$, $p_2=e_2$, $p_3=e_3$, and $p_4=(x_1,\dots,x_N)$, because M\"obius transformations are triply transitive. Choose $v_1=e_1+e_{N+1}$, $v_2=e_2+e_{N+1}$, $v_3=e_3+e_{N+1}$, and $v_4=(x_1,\dots,x_N,1)$. Then
\[
|p_1,p_2,p_3,p_4|=\sqrt{\frac{1-x_3}{1-x_2}}
\]
and
\[
|\ell_1,\ell_2,\ell_3,\ell_4|=\frac{1-x_3}{1-x_2}.
\]

To establish Theorem~\ref{T: main2} we prove the following reformulation  of Theorem~\ref{T: main2} in terms of the hyperboloid model of hyperbolic space and Lorentz transformations.

\begin{proposition}\label{P: cross-ratio}
Given two collections of  light-like lines $\{\ell_\alpha\,:\,\alpha\in \mathcal{A}\}$ and  $\{\ell'_\alpha\,:\,\alpha\in \mathcal{A}\}$, there is a positive Lorentz transformation $\phi$ with $\phi(\ell_\alpha)=\ell'_\alpha$ for each $\alpha$ in $\mathcal{A}$ if and only if $|\ell_\alpha,\ell_\beta,\ell_\gamma,\ell_\delta|=|\ell'_\alpha,\ell'_\beta,\ell'_\gamma,\ell'_\delta|$ for all ordered 4-tuples $(\alpha,\beta,\gamma,\delta)$ of distinct indices in $\mathcal{A}$.
\end{proposition}
\begin{proof}
The proposition is true when $\mathcal{A}$ has fewer than four elements, by triple transitivity of M\"obius transformations. We assume, therefore, that $\mathcal{A}$ has at least four elements.  

Suppose that $|\ell_\alpha,\ell_\beta,\ell_\gamma,\ell_\delta|=|\ell'_\alpha,\ell'_\beta,\ell'_\gamma,\ell'_\delta|$ for all ordered 4-tuples $(\alpha,\beta,\gamma,\delta)$. Choose any three indices $1$, $2$, and $3$ from the collection $\mathcal{A}$. For each $i=1,2,3$, let $v_i$ be a positive vector in $\ell_i$, and let $v'_i$ be a positive vector in $\ell'_i$, chosen such that $\langle v_1,v_2\rangle=\langle v'_1,v'_2\rangle$, $\langle v_1,v_3\rangle=\langle v'_1,v'_3\rangle$, and $\langle v_2,v_3\rangle=\langle v'_2,v'_3\rangle$. This can be achieved by adjusting the $v_i$ by positive scalar multiples. For all other lines $\ell_\alpha$, choose $v_\alpha$ to be the unique positive member of $\ell_\alpha$ for which
\[
\langle v_\alpha,v_2\rangle=-\frac{\langle v_2,v_3\rangle}{\langle v_1,v_3\rangle}.
\]
Define a positive element $v'_\alpha$ of $\ell'_\alpha$ so that $\langle v'_\alpha,v'_2\rangle$ is also equal to this quantity. Then the equation
\[
|\ell_\alpha,\ell_1,\ell_2,\ell_3|=|\ell'_\alpha,\ell'_1,\ell'_2,\ell'_3|
\]
ensures that $\langle v_\alpha,v_1\rangle = \langle v'_\alpha,v'_1\rangle$. Finally , given any pair $\alpha,\beta$ in $\mathcal{A}$, the equation
\[
|\ell_\alpha,\ell_\beta,\ell_1,\ell_2|=|\ell'_\alpha,\ell'_\beta,\ell'_1,\ell'_2|
\]
reduces to  $\langle v_\alpha,v_\beta\rangle = \langle v'_\alpha,v'_\beta\rangle$. 

The subspace spanned by $\{v_\alpha\,:\,\alpha \in \mathcal{A}\}$ is time-like because
\[
\|v_1+v_2\|^2=2\langle v_1,v_2\rangle <0;
\]
therefore Proposition~\ref{P: existence} applies to yield a Lorentz transformation $\psi$ with $\psi(v_\alpha)=v'_\alpha$ for all $\alpha$. Define $\phi$ to be whichever of the maps $\psi$ or $-\psi$ is positive. Then $\phi(\ell_\alpha)=\ell'_\alpha$ for each  $\alpha$ in $\mathcal{A}$, as required. The converse implication follows by preservation of the cross-ratio under Lorentz transformations.
\end{proof}

The map $\phi$ in Proposition~\ref{P: cross-ratio} is \emph{not} unique if and only if the time-like space $V$ spanned by the $v_\alpha$ is not equal to $\mathbb{R}^{N+1}$. In other words, $\phi$ is not unique if and only if all the $v_\alpha$ lie in a time-like Euclidean plane. In terms of the points $p_\alpha$, this occurs if and only if there is an $(N-2)$-dimensional sphere in $\mathbb{S}^{N-1}$ that contains all the points $p_\alpha$.


\begin{thebibliography}{AW1}
 
\providecommand{\href}[2]{#2}

\bibitem{Be1983} A. F. Beardon, {\it The geometry of discrete groups}, Springer, New York, 1983.

\bibitem{BeMi2006} A. F. Beardon\ and\ D. Minda, Conformal automorphisms of finitely connected regions, in {\it Transcendental dynamics and complex analysis}, 37--73, Cambridge Univ. Press, Cambridge.

\bibitem{Co1995} J. B. Conway, {\it Functions of one complex variable. II}, Springer, New York, 1995.

\bibitem{Fo1951} L. R. Ford, {\it Automorphic functions}, Chelsea, 1951.

\bibitem{Ra1994} J. G. Ratcliffe, {\it Foundations of hyperbolic manifolds}, Springer, New York, 1994. 

\end{thebibliography}
\end{document}